\documentstyle[12pt,amscd]{amsart}
\newtheorem{Theorem}{Theorem}[section] 
\newtheorem{Lemma}[Theorem]{Lemma}
\newtheorem{Corollary}[Theorem]{Corollary} \newtheorem{Proposition}[Theorem]{Proposition}

\newtheorem*{Remark}{Remark}
\newtheorem*{Example}{Example}

\def\Min{\operatorname{Min}} 
 
\def\reg{\operatorname{reg}} 
\def\In{\operatorname{in}} 
 \def\gin{\operatorname{Gin}} 
\def\chara{\operatorname{char}}
   
\def\To{\longrightarrow} 
\def\syz{\operatorname{Syz}}  
\def\Ext{\operatorname{Ext}}
\def\Hom{\operatorname{Hom}}
\def\Tor{\operatorname{Tor}}
\def\geom{\operatorname{g-reg}}
\def\Proj{\operatorname{Proj}}

\def\To{\longrightarrow}
\def\m{{\frak m}}
\def\p{{\frak p}}

\def\n{{\frak n}}
\def\O{{\mathcal O}}
\def\ZZ{{\mathbb Z}}

\begin{document}

\title{Castelnuovo-Mumford regularity\\ and related invariants}
\author{Ng\^o Vi\^et Trung}
\address{Institute of Mathematics, Vien Toan Hoc, 18 Hoang Quoc Viet, 10307 Hanoi, Vietnam}
\email{nvtrung@@math.ac.vn}
\maketitle

\begin{abstract}
These notes are an introduction to some basic aspects of the Castelnuovo-Mumford regularity and 
related topics such as weak regularity, $a^*$-invariant and partial regularities. 
\end{abstract}

\setcounter{tocdepth}{1} \tableofcontents

\section*{Introduction}

Let $R = k[x_1,...,x_n]$ be a polynomial ring over a field. Let $M$ be a finitely generated graded 
$R$-module. The structure of $M$ is best understood by a minimal graded free resolution of $M$:
$$0 \To  F_s \To  \cdots \To  F_1 \To  F_0 \To  M \To 0\ .$$
Let $b_i(M)$ denote the maximal degree of the generators of $F_i$.
The Castelnuovo-Mumford regularity or regularity of $M$ is usually defined as the number
$$\reg(M) = \max\{b_i(M)-i|\ i = 0,...,s\}.$$
 
The regularity was introduced by Mumford by generalizing geometric idea of Castelnuovo \cite{Mu}. 
Originally, it was defined for a coherent sheaf by the vanishing of the sheaf cohomology. This 
approach has led to the equivalent definition 
$$\reg(M) = \max\{a_i(M)+i|\ i \ge 0\},$$ 
where $a_i(M)$ denotes the largest non-vanishing degree of the $i$th local cohomology of $M$ with 
respect to the maximal graded ideal of $R$. 

The regularity is a measure for the complexity of the structure of graded modules \cite{BaM}, 
\cite{Va}. Indeed, several important invariants of graded rings and modules can be estimated by 
means of the regularity.
 
The aim of these notes is to introduce the reader to some basic results on the Castelnuovo-Mumford 
regularity and related invariants, which may serve a later and deeper study. The setting will be 
very general. We will work with standard graded algebra over an arbitrary Noetherian ring.
The regularity is defined by means of the local cohomology modules with respect to the ideal of 
homogeneous elements of positive degree. This is meant for possible applications, in particular, 
in the study of blow-up rings of ideals. The reader is assumed to be familiar with basic concepts 
of Commutative Algebra. 

The notes are divided in six sections. 

In Section 1 we show why the above two definitions of the regularity are equivalent and discuss 
their consequences when $R$ is a polynomial ring over a field. We also prove the characterizations 
of the regularity in terms of the Tor and Ext modules. These definitions and characterizations 
clearly indicate why the regularity is a measure for the complexity of the underlying module. 

In Section 2 we present a characterization of the regularity by means of a sequence of linear 
forms, which behave like a regular sequence in larger degree. Such a sequence is called a 
filter-regular sequence. The characterization reduces the computation of the regularity to the 
computation of the largest non-vanishing degree of simple quotient modules, which have finite 
length if the base ring is an artinian local ring. The notion of filter-regular sequence will play 
a crucial role in these notes.

In Section 3 we study the notion of weak regularity which controls the vanishing of the local 
cohomology modules along a shifted degree. By definition, to compute the regularity one needs to 
check the vanishing of infinitely many components of the the local cohomology modules. 
The weak regularity reduces this computation to just only a finite number of components.
We shall see that the geometric regularity
$$\geom(M) = \max\{a_i(M)+i|\ i \ge 1\},$$
which modifies the notion of regularity in Algebraic Geometry, is a weaker version of the weak 
regularity. The geometric regularity is  easier to handle with and can be used to estimate the 
regularity.

In Section 4 we introduce the invariant
$$a^*(M) = \max\{a_i(M)|\ i \ge 0\},$$
which is a counterpart of the regularity. In fact, if $R$ is a polynomial ring in $n$ variables, 
it can be shown that
$$a^*(M) = \max\{b_i(M)-i|\ i = 0,...,s\}-n.$$
We may view the $a^*(M)$ as the real regularity and $\reg(M)$ as the shifted regularity of $M$.
The $a^*$-invariant enjoys many interesting properties. In particular, if $R$ is an algebra over a 
a local ring, $a^*(M)$ is equal to the largest non-vanishing degree of the local cohomology 
modules of $M$ with respect to the maximal graded ideal. It is well-known that these local 
cohomology modules carry more information on the structure of $M$. 

In Section 5 we introduce the partial regularities
\begin{align*}
\reg_t(M) & = \max\{a_i(M)+i|\ i \le t\},\\
a_t^*(M) & = \max\{a_i(M)|\ i \le t\}
\end{align*}
for every integer $t \ge 0$. The partial regularities provide more specific information on the 
graded structure. In fact, if $R$ is a polynomial ring in $n$ variables, it can be shown that
\begin{align*}
\reg_t(M) & = \max\{b_i(M)-i|\ i \ge n-t\},\\
a_t^*(M) & = \max\{b_i(M)|\ i \ge n-t\}-n.
\end{align*}
In this case, the partial regularities can be viewed as the regularity and the $a^*$-invariant of 
the syzygies of $M$. We shall see that the partial regularities can be computed also by means of 
filter-regular sequences. 

In Section 6 we show how one can use the technique of Gr\"obner bases to compute the regularity of 
ideals in a polynomial ring. We show that if the variables in decreasing order form a 
filter-regular sequence for the quotient ring, then the regularity does not change when passing to 
the initial ideal with respect to the reverse lexicographic order. That is for example the case of 
the generic initial ideal. To compute the regularity of a monomial ideal we introduce simpler 
invariants which are obtained from the given ideal by certain substitution $x_i = 0,1$. These 
invariants have simple combinatorial descriptions and can be computed effectively. The regularity 
is the maximum of these invariants. As a consequence, if the characteristic of the field is zero, 
the regularity of the generic initial ideal with respect to the reverse lexicographic order is 
equal to the maximal degree of the generators.
Similar results are proved for the $a^*$-invariant and the partial regularities.

\subsection*{Acknowledgement}
The author is grateful to the organizers of the International Conference on Commutative Algebra 
and Combinatorics at Harish-Chandra Research Institute, Allahabad, December 2003, for generous 
support and hospitality. Moreover, he is indebted to Irena Swanson for her reading of the first 
draft and for the suggestions which have brought this version in life.

\section{Castelnuovo-Mumford regularity} 

Let $R = k[x_1,\ldots,x_n]$ be a polynomial 
ring over a field $k$. 
Let $M$ be a finitely generated graded $R$-module.
\smallskip

By Hilbert syzygy theorem, $M$ has a minimal graded free resolution of $M$ of finite length 
$$0 \To  F_s \To  \cdots \To  F_1 \To  F_0 \To  M \To 0\ ,$$
where the maps are all of degree zero. 
Let $b_i(M)$ denote the maximum degree of the generators of $F_i$.
The {\it Castelnuovo-Mumford regularity} or {\it regularity} of $M$ is the number
$$\reg(M) := \max\{b_i(M)-i|\ i = 0,\ldots,s\}.$$

\begin{Example} {\rm We have
\begin{itemize}
\item $\reg(R) = 0$.
\item $\reg(R/fR) = \deg(f)-1$ for any homogeneous form $f \in R$.
\end{itemize}}
\end{Example}

By the above definition,  the graded pieces of $M$ behave  similarly past the degree $\reg(M)$.
For instance, if we denote by $d(M)$ the maximal degree of the generators of $M$, then 
$$d(M) = b_0(M) \le \reg(M).$$

There are several interesting characterizations of the regularity of which the following are 
especially useful. 

\begin{Proposition} \label{Ext}
\begin{align*}
\reg(M) & = \max\{t|\ \Tor_i^R(M,k)_{t-i} \neq 0\ \text{\rm for some}\ i \ge 0\},\\
& = \max\{t|\ \Ext^i_R(M,R)_{-t-i} \neq 0 \ \text{\rm for some}\ i \ge 0\}.
\end{align*}
\end{Proposition}

\begin{pf}
By the minimality of the given resolution we have $\Tor_i^R(M,k) \cong F_i\otimes k$. Therefore, 
$$b_i(M) = \max\{t|\ \Tor_i^R(M,k)_t \neq 0 \},$$
which implies the first formula for $\reg(M)$.\par
To prove the second formula we put  $r = \reg(M)$
and $F^*_i = \Hom_R(F_i,R)$. 
Since $F_i$ has no generators of degree $> r +i$,  $F^*_i$ must be zero in degree $< -r-i$. Note 
that $\Ext_R^i(M,R)$ is the homology of the dual of the resolution of $M$ at $F^*_i$. Then 
$\Ext_R^i(M,R)_t = 0$ for $t < -r-i$. 
Now let $i$ be the largest index such that $r = b_i-i$. Then $F^*_i$ has $R(r+i)$ as a summand, 
whereas $F^*_{i+1}$ has no summand of the form $R(m)$ with $m \ge r+i$. By the minimality of the 
resolution, the summand $R(r+i)$ of $F^*_i$ must map to zero in  $F^*_{i+1}$. Moreover, nothing in 
$F^*_{i-1}$ can map on to the generator of $R(r+i)$ in $F^*_i$, so it gives a nonzero class in 
$\Ext_R^i(M,R)$ of degree $-r-i$.  Thus,
$$r = \max\{t|\ \Ext_R^i(M,R)_{-t-i} \neq 0\; \text{for some}\; i \ge 0\},$$
as required.
\end{pf}

The above definition and characterization indicate clearly why the regularity can be used as a 
measure for the complexity of the structure of $M$. For more information on this topics we refer 
to \cite{BaM}, \cite{EG}, \cite{Va}.
\smallskip

We shall see that the regularity can be defined for modules over a larger class of graded algebras 
for which we don't have finite minimal free resolutions.
\smallskip
\par From now on let $R = \oplus_{n \ge 0}R_n$ be a finitely generated standard graded ring over
a noetherian commutative ring $R_0$,  where ``standard"
means $R$ is generated by elements of $R_1$. 

For any graded $R$-module $M$ we set
$$\Gamma_{R_+}(M) = \{e \in M|\  \text{$R_te = 0$ for $t \gg 0$}\}.$$
It is easy to check that $\Gamma_{R_+}(.)$ is a covariant and left exact  functor on the category 
of modules over $R$.
The {\it local cohomology} functors are the right derived functors of $\Gamma_{R_+}(.)$ 
For $i \ge 0$ we denote by $H_{R_+}^i(M)$ the $i$th local cohomology graded module of $M$. 

\begin{Remark}
{\rm $H_{R_+}^0(M) = \Gamma_{R_+}(M) = \cup_{t \ge 0}(0_M:(R_+)^t)$, which is the intersection of 
all primary components of $0_M$ whose associated prime ideals do not contain $R_+$.}
\end{Remark}

Now let $M$ be a finitely generated $R$-module. It is known that $H_{R_+}^i(M)$ is a graded module 
with $H_{R_+}^i(M)_t = 0$ for $t \gg 0$. Moreover, if $M \neq 0$, there exists an integer $s \le 
\dim M$ such that $H_{R_+}^s(M) \neq 0$ and $H_{R_+}^i(M) = 0$ for $i > s$. For more detail see 
e.g.~\cite{BrS}. \smallskip

To study the graded structure of $H_{R_+}^i(M)$ we introduce the following notation. 
\smallskip

For any graded $R$-module $H$ with $H_t = 0$ for $t \gg 0$ we set
$$a(H) := \sup\{t|\ H_t \neq 0\}$$
with the convention $a(H) = -\infty$ if $H = 0.$
This invariant can be understood as the largest non-vanishing degree of $H$.
\smallskip

If we set
$a_i(M) := a(H_{R_+}^i(M)),$
then the {\it Castelnuovo regularity} of $M$ is defined as the number
$$\reg(M) :=\ \max\{a_i(M)+i|\ i \ge 0\}.$$ From the above properties of local cohomology modules 
one can see that $\reg(M)$ is a finite number if $M \neq 0$.

\begin{Remark} {\rm
If $M_t = 0$ for $t \gg 0$, we have $H_{R_+}^i(M) = M$ and $H_{R_+}^i(M) = 0$ for $i > 0$. 
Therefore, 
$\reg(M) = a_0(M) = \max\{t|\ M_t \neq 0\}.$
In this case, $\reg(M)$ is the largest non-vanishing degree of $M$.}
\end{Remark}

If $R$ is a polynomial ring over a field, the new definition of regularity coincides with the 
previous one. This is a consequence of the following local duality (see e.g.~\cite{BrS}).

\begin{Theorem} \label{dual}  {\rm (Local duality)}
Let $R$ be a polynomial ring over $k$ in $n$ variables. Then 
$$H_{R_+}^i(M)_t \cong \Ext_R^{n-i}(M,R)_{-t-n}$$ 
for all $i$ and $t$. 
\end{Theorem}

Now, since $\Ext_R^j(M) = 0$ for $j > n$, we obtain 
\begin{align*}
\max\{a_i(M)+i|\ i > 0\} & = \max\{t|\ \Ext_R^{n-i}(M,R)_{-t-n+i} \neq 0 \ \text{\rm for some}\ i 
= 0,...,n\}\\
& = \max\{t|\ \Ext_R^j(M,R)_{-t-j} \neq 0 \ \text{\rm for some}\ j \ge 0\}.
\end{align*}
By Proposition \ref{Ext}, this shows that the new definition of regularity coincides with the one 
in the case $R$ is a polynomial ring over a field.

The new definition of regularity has many advantages.
First, it holds for a larger class of graded algebras where not all modules have  finite free 
resolutions. Second, it gives us some flexibility in choosing the base ring.

First, if $R$ is a quotient of a graded algebra $S$, $H_{R_+}^i(M) = H_{S_+}^i(M)$ for all $i$. 
Therefore, the regularity of a finitely graded $R$-module $M$ does not change if we consider $M$ 
as a graded $S$-module.

Second, if there exists a standard graded subalgebra $A$ of $R$ such that $R$ is a finite 
$A$-module, every element of $R_+$ is integral over $A$. Therefore,
$R_+$  is contained in the radical of $A_+R$. This implies $H_{A_+}^i(M) = H_{R_+}^i(M)$ for all 
$i$. Hence $\reg(M)$ does not change if we consider $M$ as a graded $A$-module.

\begin{Remark} {\rm 
If $R$ is a standard graded algebra over a field, we may represent $R$ as a quotient ring of a 
polynomial ring over a field.
In this case, we can use the definition of regularity by means of a minimal finite resolution 
again. An alternative way to pass to the case of a polynomial ring is by taking a Noether 
normalization of $R$.}
\end{Remark}

Using the definition by means of local cohomology we can easily deduce some basic properties of 
the regularity. For instance, the regularity behaves well in a short exact sequence.

\begin{Lemma}  \label{exact}
Let $0 \To E  \To M \To F \To 0$ be an exact sequence of finitely generated $R$-modules. Then \par
{\rm (a) } $\reg(E) \le \max\{\reg(M),\reg(F)+1\}$,\par
{\rm (b) } $\reg(M) \le \max\{\reg(E), \reg(F)\},$ \par
{\rm (c) } $\reg(F) \le  \max\{\reg(M), \reg(E)-1\}.$
\end{Lemma}

\begin{pf}
We only need to consider the derived long exact sequence of local cohomology modules:
$$ \cdots \To H_{R_+}^{i-1}(F)_t \to H_{R_+}^i(E)_t \to H_{R_+}^i(M)_t \to H_{R_+}^i(F)_t \to 
H_{R_+}^{i+1}(E)_i \to \cdots\ .$$ From this it follows that
\begin{align*}
a_i(E) & \le \max\{a_{i-1}(F),a_i(M)\},\\
a_i(M) & \le \max\{a_i(E),a_i(F)\},\\
a_i(F) & \le \max\{a_i(M),a_{i+1}(E)\}.
\end{align*}
The assertion now follows from the definition of the regularity.
\end{pf}

For more basic information on the regularity we refer to \cite{BrS}, \cite{E}, \cite{EG}, 
\cite{O}.

\section{Filter-regular sequence}

In this section we present a characterization of the regularity by means of a sequence of linear 
forms, which does not involve local cohomology explicitly. This characterization will play a 
crucial role in the rest of these notes.
\smallskip

Let $R$ be a finitely generated standard graded algebra over a noetherian ring. Let $M$ be a 
finitely generated graded $R$-module.
We call a homogeneous element $z \in R$  an $M$-{\it filter-regular element} if
$(0_M:z)_t = 0$ for $t \gg 0$. This means that $z$ behaves like an $M$-regular element in higher 
degrees. If $M = R$, we delete the prefix $M$. \smallskip

It is easy to see that $z$ is an $M$-filter-regular element if and only if $z \not\in \p$ for any 
associated prime ideal $\p \not\supseteq R_+$ of $M$. Therefore, $M$-filter-regular elements of 
any degree do exist if $R_0$ is a local ring with infinite residue field. 

Using an $M$-filter-regular element one can easily compute $a_0(M)$.

\begin{Lemma} \label{a_0}
Let $z$ be any $M$-filter-regular linear form. Then
$$a_0(M) = a(0_M:z) = a(0_M:R_+).$$
\end{Lemma}

\begin{pf}
By the above characterization of filter-regular element, $\cup_{i \ge 1}(0_M:z^i)$ is the 
intersection of all primary components of $0_M$ whose associated prime ideals do not contain 
$R_+$. From this it follows that
$$\cup_{i \ge 1}(0_M:z^i) = \cup_{i \ge 1}(0_M:(R_+)^i).$$
Therefore,
$0_M:R_+ \subseteq 0_M:z \subseteq \cup_{i \ge 1}(0_M:(R_+)^i) =  H_{R_+}^0(M).$ This implies
$$a(0_M:R_+) \le a(0_M:z) \le a_0(M).$$
On the other hand, every element of the largest non-vanishing degree of 
$H_{R_+}^0(M)$ is contained in $0_M:R_+$. Therefore, $a(0_M:R_+) = a_0(M)$, which together with 
the above inequalities imply the assertion.
\end{pf}

There is the following relationship between  $\reg(M)$ and $\reg(M/zM)$.

\begin{Lemma} \label{hyperplane-1}
Let $z$ be any $M$-filter-regular linear form. Then
$$\reg(M) = \max\{a_0(M), \reg(M/zM)\}.$$
\end{Lemma}

\begin{pf}
Since $(0_M:z)_t = 0$ for $t \gg 0$, we have $H_{R_+}^i(0_M:z) = 0$ for $i \ge 1$. Now, from the 
short exact sequence $$0 \to (0_M:z) \to M \to M/(0_M:z) \to 0$$ 
we obtain $H_{R_+}^i(M) \cong H_{R_+}^i(M)$ for $i \ge 1$. From the short exact sequence 
$$0 \to M/(0_M:z) \overset z \To M \to M/zM \to 0$$
we obtain the derived long sequences of local cohomology modules:
$$\cdots \to H_{R_+}^{i-1}(M)_t \to H_{R_+}^{i-1}(M)_t \to H_{R_+}^i(M)_{t-1}\to H_{R_+}^i(M)_t 
\to H_{R_+}^i(M/zM)\to \cdots\ .$$
As a consequence, the map $H_{R_+}^i(M)_{t-1}\to H_{R_+}^i(M)_t$ is injective for $t > 
a_{i-1}(M/zM)$. Since $H_{R_+}^i(M)_t = 0$ for $t \gg 0$, this implies $H_{R_+}^i(M)_{t-1} = 0$ 
for $t > a_{i-1}(M/zM)$. 
Hence $a_i(M) \le a_{i-1}(M/zM)-1$.
The above exact sequence of local cohomology modules also implies
$a_{i-1}(M/zM) \le \max\{a_{i-1}(M), a_i(M)+1\}.$
Therefore,
$$a_i(M)+i \le a_{i-1}(M/zM) + (i-1) \le \max\{a_{i-1}(M)+ (i-1), a_i(M)+i\}$$
for $i > 0$. From this it follows that
$$\max\{a_i(M)|\ i \ge 1\} \le \reg(M/zM) \le \reg(M),$$
which immediately implies $\reg(M) = \max\{a_0(M), \reg(M/zM)\}.$
\end{pf}

Since every $M$-regular element is $M$-filter-regular, Lemma \ref{hyperplane-1} implies the 
following property. 

\begin{Corollary} 
Assume that $z$ is an $M$-regular linear form. Then 
$$\reg(M) = \reg(M/zM).$$
\end{Corollary}

\begin{pf}
The existence of an $M$-filter-regular form implies that $0_M$ has no associated prime ideals 
which contain $R_+$. This means $H_{R_+}^0(M) = 0$. Hence the assertion follows from Lemma 
\ref{hyperplane-1}.
\end{pf}

We call a sequence of homogeneous elements $z_1,...,z_s$ an $M$-{\it filter-regular sequence} if 
$z_{i+1}$ is an $M/Q_iM$-filter-regular element for $i = 0,...,s-1$, where $Q_0 = 0$ and $Q_i = 
(z_1,...,z_i)$. If $M = R$, we delete the prefix $M$. 
\smallskip

The above condition means $a((Q_iM:z_{i+1})/Q_iM) < \infty$
for $i = 0,\ldots,s-1.$ For simplicity we set
$$a(z_1,...,z_s;M) = \max\{a((Q_iM:z_{i+1})/Q_iM)|\ i = 0,...,s-1\}.$$

\begin{Remark}  {\rm Filter-regular sequences are not permutable in general. For example, let $R = 
k[x,y,z]/(x) \cap (x^2,y)$. Then $z,y$ is an $R$-filter-regular sequence, while $y,z$ is not.}
\end{Remark}

We call a homogeneous ideal $Q \subseteq R_+$ an {\it $M$-reduction} of $R_+$ if $Q$ is generated 
by linear forms and $(M/QM)_t = 0$  for $t \gg 0$. \smallskip

The above lemmas lead us to the following non-cohomological characterizations of the regularity.

\begin{Theorem} \label{regularity}
Let $z_1,...,z_s$ be an $M$-filter-regular sequence of linear forms such that $Q = (z_1,...,z_s)$ 
is an $M$-reduction of $R_+$. Then
\begin{align*}
\reg(M) & = \max\{a(z_1,...,z_s;M), a(M/QM)\}\\
& = \max\{a((Q_iM:R_+)/Q_iM)|\ i = 0,...,s-1\}.
\end{align*}
\end{Theorem}

\begin{pf}
Applying Lemma \ref{hyperplane-1} to the quotient modules
$M/Q_{i-1}M$, $i = 1,...,r$, we obtain
$$\reg(M) = \max\{a_0(M/Q_0M),...,a_0(M/Q_{s-1}M), \reg(M/QM)\}.$$
By Lemma \ref{a_0}, we have
\begin{align*}
a_0(M/Q_iM) & = a((Q_iM:z_{i+1})/Q_iM) = a((Q_iM:R_+)/Q_iM),\\
\reg(M/QM) & = a_0(M/QM) = a(QM:R_+/QM).
\end{align*} 
Therefore, the assertion is immediate.
\end{pf}

\begin{Example} 
{\rm Let $R = R_0[x_1,...,x_n]$ be a polynomial ring over $R_0$. Then $x_1,...,x_n$ is a regular 
sequence in $R$. Hence $a(x_1,...,x_n) = 0$, which implies $\reg(R) = a(R_0) = 0$.}
\end{Example}

If $R$ is a graded algebra over a local ring $(R_0,\n)$, we may assume that $R_0$ has infinite 
residue field by passing to the algebra $R\otimes R_0[u]_{\n R_0[u]}$, where $u$ is an  
indeterminate. In this case, any $M$-reduction of $R_+$ can be generated by an $M$-filter-regular 
sequence of linear forms. This is an easy consequence of the characterization of 
$M$-filter-regular elements by means of the associated primes $\p \not\supseteq R_+$ of $0_M$. 

In general, any $M$-reduction of $R_+$ can be generated by
an $M$-filter-regular sequence in a flat extension of $R$.

\begin{Lemma} \label{flat}
Let $Q$ be an $M$-reduction generated by the linear forms
$x_1,\ldots,x_s$. For $i = 1,\ldots,s$ put $z_i = 
\sum_{j=1}^su_{ij}x_j$,
where $U = \{u_{ij}|\ i, j = 1,\ldots,s \}$ is a matrix of 
indeterminates.
Put
$$R' = R[U,\det(U)^{-1}],\ R' = R \otimes_AA',\ M' = M \otimes_AA'.$$
If we view $R'$ as a standard graded algebra over $R'$ and $M'$ as a
graded $R'$-module, then $z_1,\ldots,z_s$ is an $M'$-filter-regular
sequence.
\end{Lemma}

\begin{pf}
It suffices to show that $z_1$ is $M$-filter-regular.
This will imply that $z_i$ is $(M/(z_1,...,z_{i-1})M)$-filter-regular for $i = 2,...,s$.
We will show that $z_1 \not\in P$ for any
associated prime ideal $P \not\supseteq R_+'$ of $M'$. By the definition
of $R'$, such a prime ideal $P$ must have the form  $\p R'$ for some
associated prime ideal $\p \not\supseteq R_+$ of $M$. If $Q \subseteq
\p$, then $(M/\p M)_n = 0$ for $n \gg 0$ because $M/\p M$ is
a quotient module of $M/QM$. From this it follows that there is a
number $t$ such that $(R_+)^tM \subseteq \p M$. Since
$\operatorname{ann}(M)\subseteq \p$, this implies $R_+
\subseteq \p$, a contradiction. So we
get $Q \not\subseteq \p$. Since $Q = (x_1,\ldots,x_s)$, this
implies $z_1 = u_{11}x_1 + \cdots + u_{1s}x_s \not\in \p R' = P$, as desired.
\end{pf}

Lemma \ref{flat} allows us to use Theorem \ref{regularity} for the
computation of $\reg(M)$. Indeed, since $R' = $ is a flat
extension of $R$, we have $H_{R_+'}^i(M')_n \cong
H_{R_+}^i(M)_n\otimes_AA'$ for all $n$ and $i \ge 0$,
whence
$$\reg(M) = \reg(M').$$

Let $d(M)$ denote the maximal degree of the homogeneous minimal generators of $M$. The following 
consequence of Theorem \ref{regularity}
is already known in the case $R$ is a polynomial ring.

\begin{Corollary} 
$d(M) \le \reg(M).$
\end{Corollary} 

\begin{pf}
It is easily seen that
$$d(M) = d(M/R_+M) = a(M/R_+M).$$
By Lemma \ref{flat} we may assume that $R_+$ is generated by an $M$-filter-regular sequence of 
linear forms. Applying Theorem \ref{regularity} we obtain
$$a(M/R_+M) \le \reg(M).$$
Therefore, the assertion is immediate.
\end{pf}

The characterization of the regularity by means of filter-regular sequence has been used 
successfully in the study of the reduction number and the structure of blow-up algebras 
\cite{Tr1}, \cite{Tr2}, \cite{Tr3}, \cite{Tr4}, \cite{Tr5}.

\section{Weak regularity}

By definition, to estimate the regularity we need to check the vanishing of infinitely many 
components of the local cohomology modules. We shall see that it suffices to check only a few.  
For this purpose we introduce the following notions. \smallskip

Let $R$ be a finitely generated standard graded algebra over a noetherian ring. Let $M$ be a 
finitely generated graded $R$-module. For any integer $t$ we say that $M$ is {\it $t$-regular} if 
$H_\m^i(M)_{n-i+1}=0$ for
all $i\ge 0$ and $n \ge t$.  
It is easy to check that
$$\reg(M) = \min\{t \in \ZZ|\ \text{$M$ is $t$-regular}\}.$$

We say that $M$ is {\it weakly  $t$-regular} if $H_\m^i(M)_{t-i+1}=0$ for $i \ge 0$. 
In general, weak $t$-regularity does not imply $t$-regularity.

\begin{Example}
{\rm Let $M = R = k$ be a field. Then $H_{R_+}^{i}(M) = 0$ and $H_{R_+}^{i}(E) = 0$ for $i \geq 
1$. Hence $M$ is weakly $t$-regular but not $t$-regular for $t \le -2$.}
\end{Example}

However, the two notions coincide under some mild restriction.

\begin{Theorem} \label{weak}
Let $t \ge d(M)$. Then $M$ is $t$-regular if $M$ is weakly $t$-regular.
\end{Theorem}

\begin{pf}
Using Lemma \ref{flat} we may assume that there exists an $M$-filter-regular sequence of linear 
forms $z_1,...,z_s$ such that $R_+ = (z_1,...,z_s)$.

If $s = 0$, $R_+ = 0$. In this case, $H_{R_+}^0(M) = M$ and $H_{R_+}^i(M) = 0$ for $i > 0$.
Therefore, $M$ is $t$-regular if $M_n = 0$ for $n > t$.
The latter condition is clearly satisfied because $t \ge d(M)$ and $R_+ =0$. 

If $s > 0$, we set $z = z_1$. 
As in the proof of Lemma \ref{hyperplane-1} we have the exact sequence 
$$H_{R_+}^i(M)_{n-1} \to H_{R_+}^i(M)_n \to H_{R_+}^i(M/zM)_n \to H_{R_+}^{i+1}(M)_{n-1}$$
for $i \ge 0$. Since $H_{R_+}^i(M)_{t-i+1} = 0$ and $H_{R_+}^{i+1}(M)_{t-i} = 0$ we get 
$H_{R_+}^i(M/zM)_{t-i+1} = 0$ for $i \ge 0$. Hence
$M/zM$ is weakly $t$-regular.
Since $d(M/zM) = d(M)$, using induction we may assume that $M/zM$ is $t$-regular.
This implies
$\reg(M/zM) \le t.$

As a consequence, $H_{R_+}^0(M/zM)_n = 0$ and hence 
$H_{R_+}^0(M/(0_M:z))_{n-1} =  H_{R_+}^0(M)_n$ for $n \ge t+1$.
On the other hand, since $H_{R_+}^1(0_M:z) = 0$, from the exact sequence 
$$0 \to (0_M:z) \to M \to M/(0_M:z) \to 0$$
we can derive that the map $H_{R_+}^0(M) \to  H_{R_+}^0(M/(0_M:z))$ is surjective. 
Therefore, there is a surjective map $H_{R_+}^0(M)_{n-1} \to  H_{R_+}^0(M)_n$ for all $n \ge t+1$.
Since $H_{R_+}^0(M)_{t+1} = 0$, this implies $H_{R_+}^0(M)_n = 0$ for $n \ge t+1$  or, 
equivalently, $a_0(M) \le t$.

Summing up, we get 
$$\reg(M) = \max\{a_0(M), \reg(M/zM)\} \le t.$$
Hence $M$ is $t$-regular.
\end{pf}

\begin{Corollary}
$\reg(M) = \min\{t \ge d(M)|\ \text{\rm $M$ is weakly $t$-regular}\}.$
\end{Corollary}

\begin{pf}
Since $\reg(M) \ge d(M)$, we get
\begin{align*}
\reg(M) & =  \min\{t \ge d(M)|\ \text{\rm $M$ is $t$-regular}\}\\
& = \min\{t \ge d(M)|\ \text{\rm $M$ is weakly $t$-regular}\},
\end{align*}
as desired.
\end{pf}

The regularity in algebraic geometry is defined a bit differently than in the algebraic 
case.\smallskip

For a moment let $R$ be a polynomial ring over a field and $M$ a finitely generated graded 
$R$-module. Let $X = \Proj R$. We denote by $\tilde M$ the coherent sheaf associated with $M$ on 
$X$ and by $H^i(X,\tilde M(t))$ the $i$th sheaf cohomology of the twisted $\O_X$-module $\tilde 
M(t)$, $t \in \ZZ$. Then $\tilde M$ is called $t$-{\it regular} if the map $M_n \to H^0(X, \tilde 
M(n)) $ is surjective and 
$H^i(X, \tilde M(n-i)) = 0$ for all $n > t$ and $i \ge 1$. \smallskip

There is the following relationship between the sheaf cohomology of $\tilde M$ and the local 
cohomology of $M$ (see e.g.~[BrS]).

\begin{Theorem} {\rm (Serre-Grothendieck correspondence)}
There are the exact sequence
$$ 0 \to H^0_{R_+}(M)_n \to M_t \to H^0(X, \tilde M(n)) \to H^1_{R_+}(M)_n \to 0 $$
and the isomorphisms 
$H^i(X, \tilde M(n)) \cong H^{i+1}_{R_+}(M)_n$ for $i > 0.$ 
\end{Theorem}

As a consequence, $\tilde M$ is $t$-regular if and only if $H_{R_+}^i(M)_{n-i+1} = 0$ for all $n 
\ge t$ and $i \ge 1$. This leads to the following definition. \smallskip

Let $R$ now be a finitely generated standard graded algebra over a Noetherian ring and $M$ a 
finitely generated graded $R$-module. We say that  $M$ is {\it geometrically $t$-regular} if
$H_{R_+}^i(M)_{n-i+1} = 0$ for all $n \ge t$ and $i \ge 1$.
The following result of Mumford shows that we only need to check this condition for $n = t$. 

\begin{Theorem} \label{Mumford-1}
$M$ is geometrically $t$-regular if $H_{R_+}^i(M)_{t-i+1} = 0$ for $i \ge 1$.
\end{Theorem} 

\begin{pf}
Using Lemma \ref{flat} we may assume that there exists an $M$-filter-regular sequence of linear 
forms $z_1,...,z_s$ such that $R_+ = (z_1,...,z_s)$. 
As in the proof of Theorem \ref{weak} the case $s = 0$ is trivial and if $s > 0$, we may assume 
that $H_{R_+}^i(M/zM)_{n-i+1} = 0$ for $n \ge t$ and $i \ge 1$, where $z = z_1$. 
This implies that the map
$H_{R_+}^i(M)_{n-i} \to H_{R_+}^i(M)_{n-i+1}$ is surjective for $n \ge t$ and $i \ge 1$. 
Since $ H_{R_+}^i(M)_{t-i+1} = 0$, this implies $ H_{R_+}^i(M)_{n-i+1} = 0$ for all $n \ge t$ and 
$i \ge 1$, as desired.
\end{pf}

Theorem \ref{Mumford-1} shows that geometric $t$-regularity is weaker than weak $t$-regularity and 
$t$-regularity. However, if $H_{R_+}^0(M) = 0$, the three notions coincide.
In this case, we may reformulate Theorem \ref{Mumford-1} as follows.

\begin{Corollary}
Assume that $H_{R_+}^0(M) = 0$. Then $M$ is $t$-regular if $M$ is weakly $t$-regular.
\end{Corollary}

We call the number
$$\geom(M) := \max\{a_i(M)+i|\ i \ge 1\},$$
the {\it geometric regularity} of $M$. It is clear that 
$$\geom(M) = \min\{t|\ \text{$M$ is geometrically $t$-regular}\}.$$

The regularity is related to the geometric regularity by the formula
$$\reg(M) = \max\{a_0(M),\geom(M)\}.$$
 
Comparing with the regularity, the geometric regularity has the advantage that it can be estimated 
in terms of the geometric regularity of a generic ``hyperplane plane section". For that we need 
the following observation.

Let $R$ be an algebra over an artinian local ring $R_0$. For all $t \in \ZZ$, the graded piece 
$M_t$ is an $R_0$-module of finite length, and  we can consider the Hilbert function 
$h_M(t) := \ell(M_t).$ It is well-known that $h_M(t)$ is equal to a polynomial $p_M(t)$ for $t \gg 
0$. One calls $p_M(t)$ the Hilbert polynomial of $M$. The difference $h_M(t) - p_M(t)$ can be 
expressed in terms of local cohomology modules, due to a result of Serre (see e.g. \cite{BrS}).

\begin{Theorem} \label{Serre} {\rm (Serre formula)}
Let $R$ be an algebra over an artinian local ring. Then
$$h_M(t) - p_M(t) = \sum_{i=0}^{\dim M}\ell(H_{R_+}(M)_t).$$
\end{Theorem}

This formula will be used in the proof of the following estimate for the geometric regularity, 
which is based on an idea of Mumford in  \cite[pp. 101, proof of
Theorem]{Mu}.

\begin{Theorem} \label{Mumford-2}
Let $R$ be an algebra over an artinian local ring. Let $z$ be an $M$-filter-regular linear form. 
Let $t \ge d(M)$ such that
$M/zM$ is  geometrically $t$-regular.  
Then  $M$ is geometrically $(t + p_M(t) - h_{M/L}(t))$-regular, where $L = \Gamma_{R_+}(M)$.
\end{Theorem}

\begin{pf}
We have to show that
$$\geom(M) \le t + p_M(t) - h_{M/L}(t).$$

Consider the quotient module $M/L$. Since $L_n = 0$ for $n \gg 0$, we have 
$H_{R_+}^i(L) = 0$ and hence $H_{R_+}^i(M) = H_{R_+}^i(M/L)$
for $i \ge 1$. Therefore,  $\geom(M/L) = \geom(M)$. Moreover, we have
$h_M(n) = h_{M/L}(n)$ for $n \gg 0$, which implies $p_M(t) = p_{M/L}(t)$.
On the other hand, the relation $d(M/L) \le d(M)$ shows that $t \ge d(M/L)$. Since $((L+zM)/zM)_n 
= 0$ for $n \gg 0$,  $H_{R_+}^i((L+zM)/zM) = 0$ and hence $H_{R_+}^i(M/(L+zM)) \cong 
H_{R_+}^i(M/zM)$ for $i \ge 1$. Thus, $M/(L+zM)$ is geometrically $t$-regular like $M/zM$. 
So we may replace $M$ by $M/L$. That means we may assume that $L = 0$ and $z$ is $M$-regular.

As we have seen in the proof of Lemma \ref{hyperplane-1}, 
$$a_i(M)+i \le a_{i-1}(M/zM)+i-1$$
for $i \ge 1$. Since $M/zM$ is geometrically $t$-regular,  $a_{i-1}(M/zM)+i-1\le t$ for $i \ge 2$. 
Therefore, $a_i(M)+i \le t$ for $i \ge 2$. 

Let $s$ be the smallest integer $\ge t$ such that $H_{R_+}^0(M/zM)_s = 0$. Then $M/zM$ is weakly 
$s$-regular. By Theorem \ref{weak}, this implies $\reg(M/zM) \le s$. Thus, $a_1(M)+1 \le a_0(M/zM) 
\le s$.
Since we already have $a_{i-1}(M/zM)+i-1\le s$ for $i \ge 2$, we get $\geom(M) \leq s$. 

Now, we will show that $s \le t + p_M(t) - h_M(t)$. Since $ H_{R_+}^0(M) = L = 0$ and 
$H_{R_+}^1(M/zM)_n = 0$ for $n \le t$, from the exact sequence 
$$0 \to M \overset{z} \To M \to M/zM \to 0$$
we obtain the exact sequence
$$ 0 \to H_{R+}^{0}(M/zM)_n \to  H_{R+}^{1}(M)_{n-1} \to H_{R+}^{1}(M)_n \to 0$$ 
for $n \ge t$. Since $H_{R+}^{0}(M/zM)_n \neq 0$ for $t \le n < s$, we have
$\ell(H_{R+}^{1}(M)_{n-1}) > \ell(H_{R+}^{1}(M)_n)$ for $t \le n <s$. Therefore,
$s - t \le \ell(H_{R+}^{1}(M)_t).$
On the other hand, $p_M(t) - h_M(t) = \ell(H_{R+}^{1}(M)_t)$
by  Lemma \ref{Serre}. So we obtain
$s \leq t + p_M(t) - h_M(t)$, as required.
\end{pf}

The above estimate is especially useful in finding bounds for the regularity by means of the 
degree \cite{Kl}, \cite{Li}, \cite{RTV1}, \cite{RTV2}.

\section{$a^*$-invariant} 

Let $R$ be a finitely generated standard graded algebra over a noetherian ring  and $M$ a finitely 
generated graded $R$-module. We call the number
$$a^*(M) := \max\{a_i(M)|\ i \ge 0\}$$
the $a^*$-{\it invariant} of $M$. 
In some sense, it can be considered as the ``real" regularity, whereas the Castelnuovo-Mumford 
regularity is the ``shifted" regularity of $M$. 
It is obvious that
$$a^*(M) \le \reg(M) \le a^*(M) + \dim(M).$$

The $a^*$-invariant has many interesting properties. 
First of all, if $R$ is a polynomial ring over a field, $a^*(M)$ can be characterized by means of 
the shifts in the minimal free resolution of $M$.\smallskip

\begin{Theorem} 
Let $R$ be a polynomial ring over $k$ in $n$ variables. Then
$$a^*(M) = \max\{b_i(M)|\ i \ge 0\} - n.$$
\end{Theorem}

\begin{pf}
By Theorem \ref{dual}  we have 
$$a_i(M) = \max\{t|\ \Ext_R^{n-i}(M,R)_{-t-n} \neq 0\}.$$
Therefore,
$$
a^*(M) = \max\{t|\ \Ext_R^{n-i}(M,R)_{-t-n} \neq 0\; \text{for some}\ i\ge 0\}.$$
Put 
$m = \max\{b_i|\ i \ge 0\}-n.$ For every $i \ge 0$, $F_i$ has no generators of degree $\ge m 
+n+1$, so $F^*_i = \Hom_R(F_i,R)$ must be zero in degree $\le -m-n-1$. Since $\Ext_R^i(M,R)$ is 
the homology of the dual of the resolution of $M$ at $F^*_i$,  $\Ext_R^i(M,R)_r = 0$ for $r \ge 
-m-n-1$. Now let $i$ be the largest index such that $b_i-n = m$. Then $F^*_i$ has $R(m+n)$ as a 
summand, whereas $F^*_{i+1}$ has no summand of the form $R(r)$ with $r \ge m+n$. By the minimality 
of the resolution, the summand $R(m+n)$ of $F^*_i$ must map to zero in  $F^*_{i+1}$. Moreover, 
nothing in $F^*_{i-1}$ can map on to the generator of $R(m+n)$ in $F^*_i$, so it gives a nonzero 
class in $\Ext_R^i(M,R)$ of degree $-m-n$.  Thus,
$$\max\{t|\ \Ext_R^i(M,R)_{-t-n} \neq 0\ \text{for some}\ i\ge 0\} = m,$$
as desired. 
\end{pf}

The $a^*$-invariant also controls the place where the Hilbert function $h_M(t)$ coincides with the 
Hilbert polynomials $p_M(t)$ of $M$.

\begin{Proposition} 
Let $R$ be graded algebra over an artinian local ring. Then
$h_M(n) = p_M(n)$ for $n > a^*(M).$
\end{Proposition}

\begin{pf}
This is an immediate consequence of Theorem \ref{Serre}.
\end{pf}

Similarly as for $\reg(M)$, we can compute $a^*(M)$ by means of a filter-regular sequence of 
linear forms. For this we need the following reduction.

\begin{Lemma} \label{hyperplane-2}
Let $z$ be an $M$-filter-regular linear form. Then
$$a^*(M) = \max\{a_0(M), a^*(M/zM)-1\}.$$
\end{Lemma}

\begin{pf}
We have seen in the proof of Lemma \ref{hyperplane-1} that
$$a_i(M) \le a_{i-1}(M/zM)-1 \le \max\{a_{i-1}(M)-1, a_i(M)\}$$
for $i \ge 1$. From this it follows that
$$\max\{a_i(M)|\ i \ge 1\} \le a^*(M/zM) \le a^*(M),$$
which implies the assertion.
\end{pf}

The following consequence of Lemma \ref{hyperplane-2} is sometime very useful.

\begin{Corollary} 
Assume that $z$ is an $M$-regular linear form. Then
$$a^*(M) = a^*(M/zM)-1.$$
\end{Corollary}

\begin{pf}
The existence of an $M$-filter-regular form implies $H_{R_+}^0(M) = 0$. Hence the assertion is 
immediate.
\end{pf}

For any sequence of homogeneous elements $z_1,...,z_s$ we set
$$s(z_1,...,z_r;M) = \max\{a((Q_iM:z_{i+1})/Q_iM)-i|\ i = 0,...,s-1\},$$
where $Q_i = (z_1,...,z_i)$ (so $Q_0 = 0$).

\begin{Theorem}
Let $z_1,...,z_s$ be an $M$-filter-regular sequence of linear forms such that $Q = (z_1,...,z_s)$ 
is an $M$-reduction. Then
\begin{align*}
a^*(M) & = \max\{s(z_1,...,z_s;M), a(M/QM)-s\}\\
& = \max\{a((Q_iM:R_+)/Q_iM)-i|\ i = 0,...,s-1\}.
\end{align*}
\end{Theorem}

\begin{pf}
Applying Lemma \ref{hyperplane-2} to the quotient modules
$M/Q_iM$, $i = 0,...,s-1$, we obtain
$$a^*(M) = \max\{a_0(M/Q_0M),...,a_0(M/Q_{s-1}M)-s+1, \reg(M/QM)-s\}.$$
By Lemma \ref{a_0}, we have
\begin{align*}
a_0(M/Q_iM) & = a((Q_iM:z_{i+1})/Q_iM) = a((Q_iM:R_+)/Q_iM),\\
a^*(M/QM) & = a_0(M/QM) = a(QM:R_+/QM),
\end{align*}
which implies the assertion.
\end{pf}

If $R$ is a graded algebra over a local ring, there are the local cohomology modules $H_\m^i(M)$ 
with respect to the maximal graded ideal $\m$ of $R$, which carry more information on the 
structure of $R$ than the local cohomology modules $H_{R_+}^i(M)$. For instance, $M$ is a 
Cohen-Macaulay module if $H_\m^i(M) = 0$ for $i < \dim M$. \smallskip

The vanishing of the local cohomology modules $H_{R_+}^i(M)$ and $H_\m^i(M)$ can be different. For 
instance, $H_{R_+}^d(M)$ may vanish, while we always have $H_\m^d(M) \neq 0$, where $d = \dim M$.
However, they share the same largest non-vanishing degree in the following sense. 
  
\begin{Theorem}
Let $R$ be a graded algebra over a local ring with maximal graded ideal $\m$. Then
$$\max\{a(H_\m^i(M))|\ i \ge 0\} = a^*(M).$$
\end{Theorem}

\begin{pf}
For simplicity we set 
$c(M) = \max\{a(H_\m^i(M))|\ i \ge 0\}.$

First, consider the case $M_n = 0$ for $n \gg 0$. Since $H_{R_+}(M) = M$ and $H_{R_+}^i(M) = 0$ 
for $i \ge 1$, we have $a^*(M) = a_0(M) = a(M)$.
It remains to show that $c(M) = a(M)$.
 
Set $r = a(M)$ and $m = \min\{n|\ M_n \neq 0\}$.
If $r = m$, $M$ concentrates only in degree $r$. Hence $H_\m^i(M)$ also concentrates in degree $r$ 
for all $i \ge 0$. Since $H_\m^{\dim M}(M) \neq 0$, we get $c(M) = r$.
If $r > m$, we consider the exact sequence
$$0 \longrightarrow M_r \longrightarrow M \longrightarrow M/M_r \longrightarrow 0.$$
By induction on $r-m$, we may assume that $c(M/M_r) = a(M/M_r) < r$.
Then $H_{R_+}^i(M/M_r)_n = 0$ for $n >r$ and $i \ge 0$.
Therefore, $H_\m^i(M)_n \cong H_\m^i(M_r)_n$ for $n > r$ and $i \ge 0$.
Since $c(M_r) = r$, this implies $c(M) = r$.

For the general case we may assume that $R_+$ is generated by an $M$-filter-regular sequence of 
linear forms $z_1,...,z_s$.
If $s = 0$, $R_+ = 0$ and hence $M_n = 0$ for $n \gg 0$. 
If $s > 0$, we set $z = z_1$.
We shall first show that 
$$c(M) = \max\{c(0_M:z),c(M/zM)-1\}.$$

For simplicity we set $a = c(0_M:z)$ and $b = c(M/zM)-1$.
Since $H_\m^i(0_M:z)_n = 0$ for $n > a$ and $i \ge 0$, the exact sequence 
$$0 \to (0_M:z) \to M \to M/(0_M:z) \to 0$$
implies $H_\m^i(M)_n \cong H_\m^i(M/(0_M:z))_n$ for $n > a$ and $i \ge 0$. Now, from the exact 
sequence 
$$0 \to M/(0_M:z) \overset z \To M \to M/zM \to 0$$ 
we obtain the long exact sequence  
$$\cdots \to H_\m^{i-1}(M)_{n+1} \to H_\m^{i-1}(M/zM)_{n+1} \to H_\m^i(M)_n \to H_\m^i(M)_{n+1} 
\to \cdots$$ 
for $n > a$. For $n > b$, we have $H_\m^{i-1}(M/zM)_{n+1} = 0$.
Hence the map $H_\m^i(M)_n \to H_\m^i(M)_{n+1}$ is injective for $n \ge \max\{a,b\}$.  Since 
$H_\m^i(M)_{n+1} = 0$ for $n \gg 0$, this implies $H_\m^i(M)_n = 0$ for $n > \max\{a,b\}$. 
Therefore, $c(M) \le \max\{a,b\}.$ \par

If $a < b$, we have $H_\m^{i-1}(M)_{b+1} = H_\m ^i(M)_{b+1} = 0$ for $i \ge 0$. Hence  
$H_\m^{i-1}(M/zM)_{b+1} \cong H_\m^i(M)_b$ for $i \ge 1$. Choose $i$ such that 
$H_\m^{i-1}(M/zM)_{b+1} \neq 0$. Then $H_\m^i(M)_b \neq 0$. So we obtain $c(M) \ge b$ and hence 
$c(M) = b = \max\{a,b\}$.
  
If $a \ge b$, we assume to the contrary that $c(M) < \max\{a,b\} = a$. Then $H_\m^i(M)_n = 0$ for 
$n \ge a$ and $i \ge 0$. Therefore,
$$H_\m^{i-1}(M/zM)_{a+1} = H_\m^i(M/(0_M:z))_a \cong H_\m^{i+1}(0_M:z)_a$$
for $i \ge 1$. Since $a = c(0_M:z)$, there is an index $i$ such that
$H_\m^{i+1}(0_M:z)_a \neq 0$. Therefore, 
$H_\m^{i-1}(M/zM)_{a+1} \neq 0$. This implies $b \ge a+1$, a contradiction.

Now, since $(0_M:z)_n = 0$ for $n \gg 0$, we have $c(0_M:z) = a(0_M:z) = a_0(M)$.
By induction on $s$ we may also assume that $c(M/zM) = a^*(M/zM)$. Therefore,
$$c(M) = \max\{a_0(M),a^*(M/zM)-1\}.$$
By Lemma \ref{hyperplane-2}, this implies $c(M) = a^*(M)$.
\end{pf}

The $a^*$-invariant plays an important role in the study of the Cohen-Macaulayness of graded 
algebras \cite{AHT}, \cite{Hy}, \cite{Tr2}. For more information on the $a^*$-invariant we refer 
to \cite{Sh}, \cite{Tr4}.

\section{Partial regularities}

Let $R$ be a finitely generated standard graded algebra over a noetherian ring  and $M$ a finitely 
generated graded $R$-module. Given an integer $t \ge 0$ we call the numbers
\begin{align*}
\reg_t(M) & := \max\{a_i(M)+i|\ i \le t\},\\
a_t^*(M) & := \max\{a_i(M)|\ i \le t\}.
\end{align*}
the {\it partial regularities} of $M$. 
\smallskip

The regularity and the $a^*$-invariant are special cases of the partial regularities because 
$\reg(M) = \reg_t(M),$ $ a^*(M) = a^*_t(M)$
for all $t \ge \max\{i|\ H_{R_+}^i(M) \neq 0\}$.

The meaning of the partial regularities lies on the fact that they yields specific information on 
the vanishing of the local cohomology modules and, if $R$ is a polynomial ring, on the minimal 
free resolution.

\begin{Theorem} \label{Betti}
Let $R$ be a polynomial ring over $k$ in $n$ variables. Then\par
{\rm (a) } $\reg_t(M) = \max\{b_i(M)-i|\ i \ge n-t\},$ \par
{\rm (b) } $a^*_t(M)  =  \max\{b_i(M)|\ i \ge n-t\} - n$.
\end{Theorem}

\begin{pf}
The proof is similar to the proof for the characterization of $\reg(M)$ and $a^*(M)$ in terms of 
the numbers $b_i(M)$. Hence we omit it.
\end{pf}

Theorem \ref{Betti}(a) shows that $\reg_t(M)$ is related to the following notion introduced by 
Bayer-Charalambous-Popescu \cite{BaCP}:  
$$t\text{-}\reg(M) := \max\{b_i(M)-i|\ i \ge t\}.$$
By the above theorem, we have $t\text{-}\reg(M) = \reg_{n-t}(M).$

Let $\syz_t(M)$ denote the $t$-th syzygy module of $M$ (which is defined as the kernel of the map 
$F_t \to F_{t-1}$ in the minimal free resolution of $M$ as an $R$-module). As a consequence of 
Theorem \ref{Betti} we obtain the following relationships between the partial regularity of $M$ 
and the regularity of its syzygy modules.

\begin{Corollary}  
Let $R$ be a polynomial ring over a field in $n$ variables. Then\par
{\rm (a) }  $\reg_t(M) = \reg(\syz_{n-t}(M))-n+t$,\par
{\rm (b) } $a^*_t(M) = a^*(\syz_{n-t}(M))$.
\end{Corollary}

\begin{pf}
This follows from the fact that
$$0 \to F_r \to \cdots \to F_{n-t+1} \to F_{n-t} \to \syz_{n-t}(M)$$
is a minimal free resolution of $\syz_{n-t}(M)$. 
Indeed, we have 
$$b_i(M)= b_{i-n+t}(\syz_{n-t}(M))$$
for $i \ge n-t.$
Putting $j = i-n+t$ we get 
\begin{align*} 
\reg_t(M) & = \max\{b_j(\syz_{n-t}(M))-j +n-t))|\ j \ge 0\} =  \reg(\syz_{n-t}(M))+n-t,\\
a^*_t(M) & = \max\{b_j(\syz_{n-t}(M))|\ j \ge 0\} - n = a^*(\syz_{n-t}(M)).
\end{align*}
\end{pf}

One can also compute the partial regularities by means of filter-regular sequences. 
The main point is the following reduction.

\begin{Lemma} \label{hyperplane-3}
Let $z$ be an $M$-filter-regular linear form. For any $t \ge 1$ we have:\par
{\rm (a) } $\reg_t(M) = \max\{a_0(M), \reg_{t-1}(M/zM)\}$,\par
{\rm (b) } $a_t(M) = \max\{a_0(M), a_{t-1}(M/zM)-1\}$.
\end{Lemma}

\begin{pf}
As we have seen in the proof of Lemma \ref{hyperplane-1},
$$a_i(M)+i \le a_{i-1}(M/zM) + i-1 \le \max\{a_{i-1}(M)+i-1,a_i(M)+i\}$$
for $i \ge 1$. Taking the maximum of each term for $i \le t$ we obtain
$$\max\{a_i(M)+i|\ i = 1,...,t\} \le \reg_t(M/zM) \le \reg_t(M),$$ from which (a) immediately 
follows.

The first inequalities can be rewritten as
$$a_i(M) \le a_{i-1}(M/zM) -1 \le \max\{a_{i-1}(M)-1,a_i(M)\},$$
which implies (b) similarly as above.
\end{pf}

\begin{Corollary}
Assume that $z$ is an $M$-regular element. Then\par
{\rm (a) } $\reg_t(M) =  \reg_{t-1}(M/zM)$,\par
{\rm (b) } $a_t(M) =  a_{t-1}(M/zM)-1$.
\end{Corollary}

There are two ways for the computation of $\reg_t(M)$ and $a_t^*(M)$ by using $M$-filter-regular 
sequences of length $t$ and $t+1$.

\begin{Theorem} \label{partial-1} 
Let $z_1,...,z_t$ be an $M$-filter-regular sequence of linear forms. 
Let $Q_0 = 0$ and $Q_i = (z_1,...,z_i)$, $i = 1,...,t$. Then \par
{\rm (a) } $\reg_t(M) = \max\{a((Q_iM:R_+)/Q_iM)|\ i = 0,...,t\}.$ \par  
{\rm (b) } $a_t^*(M) = \max\{a((Q_iM:R_+)/Q_iM)-i|\ i = 0,...,t\}.$ 
\end{Theorem}

\begin{pf}
Applying Lemma \ref{hyperplane-3} to the quotient module $M/Q_iM$ we get
\begin{align*}
\reg(M) &= \max\{a_0(M/Q_iM)|\ i = 0,...,t\},\\
a^*(M)  & = \max\{a_0(M/Q_iM)-i|\ i = 0,...,t\}.
\end{align*}
Since $a_0(M/Q_iM) = a((Q_iM:R_+)/Q_iM)$ by Lemma \ref{a_0}, this implies  the assertion.
\end{pf}

\begin{Theorem} \label{partial-2}
Let $z_1,...,z_{t+1}$ be an $M$-filter-regular sequence of linear forms. Then \par
{\rm (a) } $\reg_t(M) = a(z_1,\ldots,z_{t+1};M).$ \par  
{\rm (b) } $a_t^*(M) = s(z_1,\ldots,z_{t+1};M).$ 
\end{Theorem}

\begin{pf}
The assertion can be proved similarly as in the proof of Theorem \ref{partial-1} because 
$a_0(M/Q_iM) = a((Q_iM:z_{i+1})/Q_iM)$ by Lemma \ref{a_0}. 
\end{pf}

For more information on partial regularities we refer to \cite{BaCP}, \cite{Tr6}.

\section{Regularity of ideals}

Let $R = k[x_1,...,x_n]$ be a polynomial ring over a field $k$ and $I$ a homogeneous ideal in $S$. 
We will show how one can use Gr\"obner basic technique to compute the regularity of $I$.
The method was developed by Bayer and Stillman \cite{BaS1}, \cite{BaS2}. 
It was later applied in \cite{Tr6} for the computation of the partial regularities. 

First, we observe that the computation of the regularities of $I$ can be reduced to the 
computation of $R/I$.
Indeed, from the minimal free resolution of $R/I$ we can deduce that
\begin{align*}
\reg_t(I) & = \reg_t(R/I) + 1,\\
a^*_t(I) & = a^*_t(R/I).
\end{align*}

We will use the reverse lexicographic order exclusively. 
Let $\In(I)$ denote the initial ideal of $I$.
Then there is the following relationship between $I$ and $\In(I)$.

\begin{Lemma} \label{BaS} 
Let $i = n,\ldots,1$. For every integer $m \ge 0$, 
$$[(I,x_n,\ldots,x_{i+1}):x_i]_m = (I,x_n,\ldots,x_{i+1})_m$$
if and only if  
$$[(\In(I),x_n,\ldots,x_{i+1}):x_i]_m = (\In(I),x_n,\ldots,x_{i+1})_m.$$ 
\end{Lemma}

\begin{pf}
Let $f = g + x_nh_n + \cdots + x_{i+1}h_{i+1}$ be an arbitrary polynomial in 
$(I,x_n,...,x_{i+1})$, where $g \in I$. 
If $f \not\in (x_n,...,x_{i+1})$, then $\In(f)$ does not contain the 
variables $x_n,...,x_{i+1}$. This is a property of the reverse lexicographic order. In this case, 
$\In(f) = \In(g)$. This shows that
$\In(I,x_n,...,x_{i+1})$ is contained in $(\In(I),x_n,...,x_{i+1}).$ Since the inverse containment 
is trivial,
we obtain 
$$\In(I,x_n,...,x_{i+1}) = (\In(I),x_n,...,x_{i+1}).$$
From this it follows that
$$\In((I,x_n,...,x_{i+1}):x_i) = (\In(I),x_n,...,x_{i+1}):x_i.$$
Since $\In((I,x_n,...,x_{i+1}):x_i)_m = \In(I,x_n,...,x_{i+1})_m$ if and only if 
$[(I,x_n,...,x_{i+1}):x_i]_m = (I,x_n,...,x_{i+1})_m$, this implies the assertion.
\end{pf}

The assertion of Lemma \ref{BaS} can be reformulated as
\begin{align*}
& a\big (((I,x_n,\ldots,x_{i+1}):x_i)/(I,x_n,\ldots,x_{i+1})\big) =\\  &\quad\quad\quad  
a\big(((\In(I),x_n,\ldots,x_{i+1}):x_i)/(\In(I),x_n,\ldots,x_{i+1})\big).
\end{align*}
Since this holds for $i = n,...,1$, we obtain 
\begin{align*}
a(x_n,...,x_j; R/I) & = a(x_n,...,x_j; R/\In(I)),\\
s(x_n,...,x_j; R/I) & = s(x_n,...,x_j; R/\In(I)).
\end{align*}

Now we can use the characterization of the regularities by means of filter-regular sequences to 
compare the regularities of $R/I$ and $R/\In(I)$.

\begin{Theorem}  \label{in}
Let $t = 0,...,n-1$. Assume that $x_n,\ldots,x_{n-t-1}$ is a filter-regular sequence in 
$R/\In(I)$. Then\par
{\rm (a) } $\reg_t(I) = \reg_t(\In(I))$,\par
{\rm (b) } $a_t^*(I) = a_t^*(\In(I))$.\par
\end{Theorem}

\begin{pf}
The assumption means 
$$a(x_n,...,x_{n-t-1};R/I) =a(x_n,...,x_{n-t-1};R/\In(I)) < \infty.$$
Therefore,  $x_n,\ldots,x_{n-t-1}$ is also filter-regular  in $(R/I)$.
Now we can apply Theorem \ref{partial-2} to the sequence $x_n,...,x_{n-t-1}$ in $R/I$ and 
$R/\In(I)$ and obtain
\begin{align*} 
\reg_t(R/I)& = a(x_n,...,x_{n-t-1}; R/I)\\
& = a(x_n,...,x_{n-t-1}; R/\In(I)) = \reg_t(R/\In(I)),\\ 
a_t(R/I) & = s(x_n,...,x_{n-t-1}; R/I)\\
& = s(x_n,...,x_{n-t-1}; R/\In(I)) = a_t(R/\In(I)),
\end{align*}
which imply the assertion.
\end{pf}

\begin{Corollary}  
Let $k$ be an infinite field. Let $\gin(I)$ denote the generic initial ideal of $I $ with respect 
to the reverse lexicographic order. Then \par
{\rm (a) } $\reg_t(I) = \reg_t(\gin(I))$,\par
{\rm (b) } $a^*_t(I) = a^*_t(\gin(I))$. 
\end{Corollary}

\begin{pf}
Since $k$ is infinite, we may use a generic choice of the variables to assume that 
$x_n,\ldots,x_1$ is a filter-regular sequence in $R/I$. Since $\gin(I)$ is the initial ideal of 
$I$ for generic choices of the variables, the assertion follows from Theorem \ref{in}.
\end{pf}

In the following, we will present a practical method for the computation of the regularities of 
$I$ which is based on the substitutions $x_i = 0,1$ of $\In(I)$ for some variables $x_i$.
\smallskip

Let $J$ be an arbitrary monomial ideal. For $i = 0,\ldots,n-1$ set
$R_i := k[x_1,\ldots,x_{n-i}].$
Let $J_i$ be the ideal of $R_i$ obtained from $J$ by the substitution 
$x_n = \cdots = x_{n-i+1} = 0.$
Let $\tilde J_i$ denote the ideal of $R_i$ obtained from $J_i$ by the substitution $x_{n-i} = 1.$ 
\smallskip

The ideals $J_i$ and $\tilde J_i$ can be easily computed from the generators of $J$.
In fact, if $J = (f_1,\ldots,f_s)$, where $f_1,\ldots,f_s$ are monomials in $S$, then $J_i$ is 
generated by the monomials $f_j$ not divided by any of the variables $x_{n-i+1},\ldots,x_n$ and 
$\tilde J_i$ by those monomials obtained from the latter by setting $x_{n-i} = 1$.

\begin{Lemma}  \label{c_i}
$a(\tilde J_i/J_i) = a\big(((J,x_n,\ldots,x_{n-i+1}):x_{n-i})/(J,x_n,\ldots,x_{n-i+1})\big).$
\end{Lemma}

\begin{pf}
By the definition of the ideals $J_i$ and $\tilde J_i$ we have
$$\tilde J_i/J_i \cong \cup_{j\ge 1}((J,x_n,...,x_{n-i+1}):x_{n-i}^j)/(J,x_n,...,x_{n-i+1}).$$
Therefore, it suffices to show that
\begin{align*}
& a\big(\cup_{j\ge 1}((J,x_n,...,x_{n-i+1}):x_{n-i}^j)/(J,x_n,...,x_{n-i+1})\big) =\\
& \quad\quad\quad a\big(((J,x_n,...,x_{n-i+1}):x_{n-i})/(J,x_n,...,x_{n-i+1})\big).
\end{align*}
But this follows from the fact that every element of the largest non-vanishing degree of 
$$\cup_{j\ge 1}((J,x_n,...,x_{n-i+1}):x_{n-i}^j)/(J,x_n,...,x_{n-i+1})$$
is also contained in $((J,x_n,...,x_{n-i+1}):x_{n-i})/(J,x_n,...,x_{n-i+1}).$
\end{pf}

For $i = 0,...,n$ we set
$$c_i(I)  :=  a(\tilde J_i/J_i),$$
where $J = \In(I)$. Note that 
the ideal $J_i$ can be obtained by first evaluating the ideal $I$ at the substitution $x_n = 
\cdots = x_{n-i+1} = 0$ and then computing the initial ideal of the evaluated ideal. 

There is the following characterization of the regularities of $R/I$ by means of the invariants 
$c_i(I)$.

\begin{Theorem} \label{evaluation}
Assume that $c_i(I) < \infty$ for $i = 0,...,t$. Then\par
{\rm (a) } $\reg_t(I) = \max\{c_i(I)|\ i = 0,\ldots,t\}+1,$\par
{\rm (b) } $a_t^*(I) = \max\{c_i(I)-i|\  i = 0,\ldots,t\}.$
\end{Theorem}

\begin{pf}
By Lemma \ref{c_i}, the assumption implies that $x_n,...,x_{n-t-1}$ is filter-regular in 
$R/\In(I)$. Now we can apply Theorem \ref{in} and obtain
\begin{align*}
\reg_t(I) & = \reg_t(\In(I)) = \reg_t(R/\In(I)) + 1,\\
a_t^*(I) & = a_t^*(\In(I)) = a_t^*(R/\In(I)).
\end{align*}
Therefore, the assertion follows from Theorem \ref{partial-2} and Lemma \ref{c_i}.
\end{pf}

\begin{Example}
{\rm  
Let $R = k[x_1,x_2,x_3,x_4]$ and 
$$I = (x_1x_2-x_3x_4,x_1x_3^2-x_2^3,x_1^2x_3-x_2^2x_4,x_1^3-x_2x_4^2).$$
Then $\dim R/I = 2$ and $\In(I) = (x_1x_2,x_2^3,x_1^2x_3,x_1^3)$. We have
\begin{align*}
J_0 & =  (x_1x_2,x_2^3,x_1^2x_3,x_1^3),\ \tilde J_0 = J_0,\\
J_1 & =  (x_1x_2,x_2^3,x_1^2x_3,x_1^3),\ \tilde J_1 = (x_1x_2,x_2^3,x_1^2),\\
J_2 & = (x_1x_2,x_2^3,x_1^3),\ \tilde J_2 = (x_1).
\end{align*}
From this it follows that 
$c_0(I) = -\infty,\ c_1(I) = 2,\ c_2(I) = 2,$
which implies $\reg(R/I) = \reg_2(R/I) = 2$.}
\end{Example}

We would like to mention that the numbers $c_i(I)$ can be described combinatorially in terms of 
the lattice vectors of the generators of $\In(I)$. This description together with the formulae of 
Theorem \ref{evaluation} give an effective method for the computation of the regularities of 
$R/I$. For detail we refer to \cite{BeG}, \cite{Tr7}.
\smallskip

Finally, we will use the above method to characterize the regularities of $I$ in terms of the 
minimal generators of $\gin(I)$.
\smallskip

Let $\cal B$ be the Borel subgroup of GL$(n,k)$ consisting of the upper triangular invertible 
matrices. Let $B$ act on $R$ by $\pi(f) = f(\pi\cdot (x_1,...,x_n)^T)$ for all $\pi \in B$, $f \in 
R$. A monomial ideal $I$ is called {\it Borel-fixed} if $g(I) = I$ for all $\pi \in \cal B$.  Our 
interest in this notion arises from the fact that $\gin(I)$ is a Borel-fixed ideal (see 
\cite{Ga}). \smallskip

In the case of characteristic zero, Borel-fixed ideals can be characterized as follows (see 
\cite{BaS1}). 

\begin{Proposition} \label{Borel} 
Assume that $\chara(k) = 0$. Let $J$ be a monomial ideal. Then $J$ is Borel-fixed if and only if 
whenever $x_1^{p_1}\cdots x_n^{p_n} \in J$, then 
$$x_1^{p_1}\cdots x_i^{p_i+q}\cdots x_j^{p_j-q}\cdots x_n^{p_n} \in J$$
for each $1 \le i < j \le n$ and $0 \le q \le p_j$. 
\end{Proposition}

For $A = (p_1,...,p_n)$ let $x^A$ denote the monomial $x_1^{p_1}...x_n^{p_n}$ and 
$$m(x^A) = \max\{i|\ p_i \neq 0\}.$$ 
Moreover, we denote by $\Min(J)$ the set of the minimal monomial generators of $J$.

\begin{Lemma}  \label{degree}
Assume that $\chara(k) = 0$. Let $J$ be a Borel-fixed monomial ideal. 
Then 
$$a(\tilde J_i/J_i) = \max\{\deg(x^A)|\ x^A \in \Min(J),\  m(x^A)= n-i\}-1.$$
\end{Lemma}

\begin{pf}  
Let $r  = \max\{\deg(x^A)|\ x^A \in \Min(A),\ m(x^A) = n-i\}-1$. 
Let $x^A$ be an arbitrary element of $\Min(J)$ with $m(x^A) = n-i$. Write $x^A = x^Bx_{n-i}$. Then 
$x^B \in \tilde J_i \setminus J_i$.
From this it follows that $r \le a(\tilde J_i/J_i)$.
It remains to show that $r \ge a(\tilde J_i/J_i)$.

Assume to the contrary that there is a monomial $x^C \in \tilde J_i \setminus J_i$ with $\deg x^C 
> r$.
Then there exists an integer $m \ge 1$ such that $x^Cx_{n-i}^m \in J$ and $m(x^Cx_{n-i}^m) = n-i$.
Since $\deg x^Cx_{n-i}^m \ge r+2$, $x^Cx_{n-i}^m$ is not a minimal generator of $J$. Therefore we 
can find a monomial $x^D \in J$ such that $x^Cx_{n-i} = x^Dx_h$ for some $h \le n-i$. Without 
restriction we may assume that 
$m$ is the smallest degree such that $x^Cx_{n-i}^m \in J$. 
Since $x^Cx_{n-i}^{m-1} \not\in J$, $h \neq i$. Thus, $x^D$ is divisible by $x_{n-i}^m$ and we may 
write $x^D = x^Ex_{n-i}^m$. It follows that $x^C = x^Ex_h$. By Lemma \ref{Borel}, this implies 
$x^C \in J$, a contradiction. 
\end{pf} 

\begin{Theorem} \label{gin}
Assume that $\chara(k) = 0$. Then\par
{\rm (a) } $\reg_t(I) = \max\{\deg x^A|\ x^A \in \Min(\gin(I)),\ m(x^A) \ge  n-t\}$,\par
{\rm (b) } $a_t^*(I) = \max\{\deg(x^A)+m(x^A)|\ x^A\in \Min(\gin(I)),\ m(x^A) \ge  n-t\}-n-1$. 
\end{Theorem}

\begin{pf}
Set $J = \gin(I)$. By Lemma \ref{degree},  
$$c_i(I) =  \max\{\deg(x^A)|\ x^A \in \Min(\gin(I)),\  m(x^A) = n-i\}-1 < \infty.$$
It is easy to check that
\begin{align*}
& \max\{c_i(I)|\ i = 0,...,t\} = \max\{\deg x^A|\ x^A \in \Min(\gin(I)),\ m(x^A) \ge  n-t\}-1,\\
& \max\{c_i(I)-i|\ i = 0,...,t\} =\\
&\quad\quad\quad \max\{\deg(x^A)+m(x^A)|\ x^A\in \Min(\gin(I)),\ m(x^A) \ge  n-t\}-n-1.
\end{align*}
Therefore, the assertion follows from Theorem \ref{evaluation}.
\end{pf}

\begin{Corollary} 
Assume that $\chara(k) = 0$. Then\par
{\rm (a) } $\reg(I) = \max\{\deg(x^A)|\ x^A\in \Min(\gin(I))\}$,\par
{\rm (b) } $a^*(I) = \max\{\deg(x^A)+m(x^A)|\ x^A\in \Min(\gin(I))\}-n-1$.
\end{Corollary}

\begin{pf}
Since $\reg(I) = \reg_n(I)$ and $a^*(I) = a^*_n(I)$, the assertion follows from Theorem \ref{gin}.
\end{pf}

\end{document}